\documentclass[10pt,twoside]{article}
\usepackage{graphicx}
\usepackage{amsmath}
\usepackage{Latex-document}

\def\C{{\mbox{\rm\kern.24em
\vrule width.03em height1.43ex depth-.052ex \kern-.26em C}}}
\def\QSet{\mbox{\rm\kern.24em
\vrule width.03em height1.48ex depth-.051ex \kern-.26em Q}}
\def\Z{{\mbox{\rm\kern.25em
\vrule width.03em height0.57ex depth0ex \kern.033em \vrule
width.03em height1.52ex depth-0.96ex \kern-.338em Z}}}
\def\N{{{\mbox{\rm I\kern-.2em N}}_0}}
\def\R{{\mbox{\rm I\kern-.22em R}}}

\newtheorem{theorem}{Theorem}[section]

\renewcommand{\thesection}{\arabic{section}}
\renewcommand{\thesubsection}{\thesection.\arabic{subsection}.}

\def\addsec{\addtocounter{section}{1} \setcounter{theorem}{0}}

\title{\bf  Singular Integrals Meet Modulation \vskip -1mm Invariance \vskip 6mm}

\author{C. Thiele\vspace*{-0.5cm}\thanks{Department of Mathematics, UCLA,
Los Angeles, CA 90095-1555, USA. E-mail: thiele@math.ucla.edu}}

\date{\vspace{-8mm}}

\begin{document}

\maketitle

\thispagestyle{first} \setcounter{page}{721}

\begin{abstract}

\vskip 3mm

Many concepts of Fourier analysis on Euclidean spaces rely on the
specification of a frequency point. For example classical
Littlewood Paley theory decomposes the spectrum of functions into
annuli centered at the origin. In the presence of structures which
are invariant under translation of the spectrum (modulation) these
concepts need to be refined. This was first done by L. Carleson in
his proof of almost everywhere convergence of Fourier series in
1966. The work of M. Lacey and the author in the 1990's on the
bilinear Hilbert transform, a prototype of a modulation invariant
singular integral, has revitalized the theme. It is now subject of
active research which will be surveyed in the lecture. Most of the
recent related work by the author is joint with C. Muscalu and T.
Tao.

\vskip 4.5mm


\end{abstract}

\vskip 12mm

\section*{1. Multilinear singular integrals}\label{multilinear_section} \addsec

\vskip-5mm \hspace{5mm}

A basic example for the notion of {\it  singular integral} is a convolution operator
\begin{equation}\label{convolution_operator}
Tf(x)=K*f(x)=\int K(x-y)f(y)\, dy
\end{equation}
whose convolution kernel $K$ is not absolutely integrable. If $K$
was absolutely integrable then we had trivially an a priori
estimate
\begin{equation}\label{linearp}
\|K*f\|_p\le \|K\|_1\|f\|_p
\end{equation}
for $1\le p\le \infty$. This follows by standard interpolation
techniques from the two endpoints $p=1,\infty$, which are true by
trivial manipulations.

A basic point of singular integral theory is that an estimate of the form (\ref{linearp}) may prevail for
$1<p<\infty$ with a constant $C_{p,K}$ instead of $\|K\|_1$ on the right hand side, if $K$ is not absolutely
integrable and the integral (\ref{convolution_operator}) is only defined in a distributional (principal value)
sense. The most prominent example on the real line (indeed, all operators in this article will act on functions on
the real line) is the Hilbert transform with $K(x)=1/x$.

Taking formally Fourier transforms, one can write
(\ref{convolution_operator}) as multiplier operator:
\begin{equation}\label{multiplier_operator}
\widehat{Tf}(\xi)=\widehat{K}(\xi)\widehat{f}(\xi)=:m(\xi)\widehat{f}(\xi).
\end{equation}
For the purpose of this survey a sufficiently interesting class of
singular integrals is described in terms of the multiplier $m$ by
imposing the symbol estimates
\begin{equation}\label{linear_symbol}
 (d/d\xi) ^{\alpha} m(\xi)\le C|\xi|^{-\alpha}
\end{equation}
for $\alpha=0,1,2$. We define the dual bilinear form
\begin{equation}\label{bilinear}
\Lambda(f_1,f_2)=\int (Tf_1(x)) f_2(x)\, dx=
\int_{\xi_1+\xi_2=0}\widehat{f_1}(\xi_1)\widehat{f_2}(\xi_2)m(\xi_1)\,
d\sigma
\end{equation}
where $d\sigma$ is the properly normalized Lebesgue measure on the
hyperplane $\xi_1+\xi_2=0$. The natural generalization of estimate
(\ref{linearp}) using duality of $L^p$ spaces then takes the form
\begin{equation}\label{bilinear_form}
|\Lambda(f_1,f_2)|\le C_{p_1}\|f_1\|_{p_1}\|f_2\|_{p_2}
\end{equation}
with $1/p_1+1/p_2=1$.

Estimate (\ref{bilinear_form}) can be related to square function
estimates which are fundamental in singular integral theory. Let
$(\psi_j)_{j\in \Z}$ be a family of functions such that
$m_j:=\widehat{\psi}_j$ is supported in the ball $B(0,2^j)$ of
radius $2^j$ around $0$, vanishes on $B(0,2^{j-2})$, and satisfies
the symbol estimates (\ref{linear_symbol}) uniformly in $j$. By
square function estimate we mean the inequality
\begin{equation}\label{square_function}
\| (\sum_j |f*\psi_j|^2)^{1/2}\|_{p} \le C_p \|f\|_p\ \
\end{equation}
which holds for $1<p<\infty$. Now let $m$ be any multiplier
satisfying (\ref{linear_symbol}). It is easy to split it as
$m(\xi_1)=\sum_j\widehat{\psi}_{1,j}(\xi_1)\widehat{\psi}_{2,j}(-\xi_1)$
for two families $\psi_{1,j}$ and $\psi_{2,j}$ as in the square
function estimate. Then we have
$$|\Lambda(f_1,f_2)|=|\sum_j\int (f_1*\psi_{1,j})(x) (f_2*\psi_{2,j})(x) dx|. $$
Moving the sum inside the integral and applying Cauchy-Schwarz,
H\"older, and (\ref{square_function}) we obtain
(\ref{bilinear_form}).

A natural generalization (see \cite{meyerc}) of (\ref{bilinear})
to multilinear forms is
\begin{equation}\label{multilinear_form}
\Lambda(f_1,\dots, f_n)=
\int_{\xi_1+\dots+\xi_n=0}m(\xi_1,\dots,\xi_{n-1})\prod_{j=1}^n
\widehat{f}_j(\xi_j)\, d\sigma
\end{equation}
with multipliers $m$ satisfying
\begin{equation}\label{classical_symbol}
 \partial^{\alpha} m(\xi')\le C|\xi'|^{-|\alpha|}.
\end{equation}

Here $\xi'=(\xi_1,\dots,\xi_{n-1})$ and $\alpha$ runs through all
multi- indices up to some order $N$. Note that the special role of
the index $n$ in the above is purely notational. The natural
estimates to ask for are
\begin{equation}\label{paraproduct}
|\Lambda(f_1,\dots,f_n)|\le C_{p_1,\dots,p_{n-1}} \prod_{j=1}^n
\|f_j\|_{p_j}
\end{equation}
for $1<p_j<\infty$ with $\sum_j 1/{p_j}=1$. In the special case
that $m$ is constant, $\Lambda(f_1,\dots,f_n)$ is a multiple of
the integral of the pointwise product of the functions $f_j$ and
estimate (\ref{paraproduct}) is simply H\"older's inequality.

We sketch a proof of (\ref{paraproduct}). Without destroying the
symbol estimates, we can split $m$ into a finite sum of
multipliers, each supported on a narrow cone with tip at the
origin. Thus assume $m$ is supported on such a cone consisting of
rays having small angle with a vector $\eta'$.

We may assume by symmetry that ${\eta'}_1=1$ is the maximal
component of $\eta'$. Then we can split $m$ into pieces $m_j$
satisfying (\ref{classical_symbol}) uniformly and supported in
$$(B(0,2^j)\setminus B(0,2^{j-2}))\times B(0,2^{j+n})^{n-2}. $$
Introduce $\eta_n$ such that $\sum_j\eta_j=0$. By symmetry among
the indices larger than $1$ we may assume $\eta_2\ge 1/n$. Then it
is easy to arrange (see Figure ``Cone'') the support of $m_j$ to
be in
$$(B(0,2^j)\setminus B(0,2^{j-2}))\times
(B(0,2^{j+n})\setminus B(0,2^{j-n}))\times B(0,2^{j+n})^{n-3}. $$

\setlength{\unitlength}{1.2mm}
\begin{figure}\label{cone}
\begin{center}
\begin{picture}(50,50)
\put(20,5){\line(0,1){30}}

\put(0,20){\line(1,0){40}}

\put(20,20){\vector(3,2){35}}

\put(60,40){$\eta'$}

\put(43,18){$\xi_1$}

\put(23,35){$\xi_2$}

\put(36,26){$m_j$}

\put(21.5,21){\circle{0.8}}

\put(23,22){\circle{1.6}}

\put(26,24){\circle{3.2}}

\put(32,28){\circle{6.4}}

\put(44,36){\circle{12.8}}
\end{picture}
\caption{``Cone''}
\end{center}
\end{figure}

Using smoothness of the multiplier $m_j$ we may use Fourier
expansion to write it as rapidly converging sum of multipliers of
elementary tensor form
$$\widehat{\psi}_{1,j}(\xi_1)\widehat{\psi}_{2,j}(\xi_2)\prod_{l=3}^n
\widehat{\phi}_{l,j}(\xi_l)$$ with
$\xi_n=-\sum_{j=1}^{n-1}\xi_{n-1}$. The symbol estimates prevail
for these elementary tensors, and thus we observe
\begin{equation}\label{phi}
(d/d\xi)^\alpha(\phi_{l,j})(\xi)\le C 2^{-\alpha j}
\end{equation}
for all derivatives up to order $N$. Observe that
$\widehat{\psi}_{l,j}$ are essentially as in
(\ref{square_function}), and $\widehat{\phi}_{l,j}$ are similar
but fail to be supported away form the origin. Applying the
elementary tensor multiplier form to $f_1,\dots,f_2$ is the same
as applying a constant multiplier to
$\psi_{1,j}*f_1,\dots,\phi_{n,j}*f_n$. Estimate
(\ref{paraproduct}) then follows from
$$
\sum_{j} \int \prod_{l=1}^2 ({\psi}_{l,j}*{f}_l)(x)\ \prod_{l=3}^n
{\phi}_{l,j}* {f}_l(x) \, d\sigma$$
$$\le C \prod_{l=1}^2 \|(\sum_j |f_l* \psi_{l,j}|^2)^{1/2} \|_{p_l}
\prod_{l=3}^\infty  \|\sup_j | \|f_l* \phi_{l,j}| \|_{p_l} \le C \prod_{l=1}^n \|f_l\|_{L^p_l}. $$
Here we have
used for $l=1,2$ the square function estimate (\ref{square_function}) and for $l>2$ the equally fundamental Hardy
Littlewood maximal inequality
$$\| \sup_j |f*\phi_{l,j}|\|_{L^p} \le C_p \|f\|_p\ \ $$
which is valid due to (\ref{phi}).

\section*{2. Modulation invariance}\setzero \addsec

\vskip-5mm \hspace{5mm}

Modulation $M_\eta$ with parameter $\eta\in \R$ is defined to be
multiplication by a character:
\[M_\eta f(x): =f(x)e^{2\pi i \eta x}. \]
This amounts to a translation of the Fourier transform of $f$.

We shall be interested in multilinear forms $\Lambda$ which have
modulation symmetries in the sense
\begin{equation}\label{modulation_symmetry}
\Lambda(f_1,\dots, f_n)=\Lambda(M_{\eta_1}f_1,\dots M_{\eta_n}f_n)
\end{equation}
for all vectors $\eta=(\eta_1,\dots,\eta_n)$ in a subspace
$\Gamma$ of the hyperplane given by $\sum\eta_j=0$.

If $\Lambda$ is given in multiplier form (\ref{multilinear_form}),
then (\ref{modulation_symmetry}) is equivalent to a translation
symmetry of the multiplier $m$:
\begin{equation}\label{translate_multiplier}
m(\xi_1,\dots,\xi_n)=m(\xi_1+\eta_1,\dots, \xi_n+\eta_n).
\end{equation}

Such a symmetry with nontrivial $\eta$ is inconsistent with the
symbol estimates  (\ref{classical_symbol}) unless $m$ is constant.
Namely, by iterating (\ref{translate_multiplier}), any point with
nonvanishing derivative of $m$ can by translated to a point far
away from the origin,  until the value of the derivative, which
remains constant at the translated points, contradicts
(\ref{classical_symbol}).

A natural replacement for (\ref{classical_symbol}) in the presence
of modulation symmetry along vectors in $\Gamma$ has been
introduced by Gilbert/Nahmod \cite{gilbertnahmod}:
\begin{equation}\label{new_symbol}
 \partial^{\alpha} m(\xi')\le C{\rm dist}(\xi',\Gamma')^{-|\alpha|}.
\end{equation}
Here $\Gamma'$ is the projection of $\Gamma$ onto the first $n-1$
coordinates. Figure ``Circles'' indicates the regions in which
multipliers of the form (\ref{new_symbol}) can be thought of as
being essentially constant.

\setlength{\unitlength}{1.2mm}
\begin{figure}\label{``Circles''}
\begin{center}
\begin{picture}(60,64)
\put(30,30){\line(0,1){28}}

\put(30,30){\line(-5,-3){24}}

\put(30,30){\line(5,-3){24}}

\put(30,30){\vector(1,2){15}}

\put(30,30){\line(-1,-2){11}}

\put(42,60){$\eta$}

\put(30,60){$\xi_1$} \put(3,13){$\xi_2$} \put(56,13){$\xi_3$}

\multiput(31,29.5)(0.5,1){28}{\circle{0.625}}
\multiput(31,29.5)(-0.5,-1){21}{\circle{0.625}}

\multiput(32,29)(1,2){15}{\circle{1}}
\multiput(32,29)(-1,-2){10}{\circle{1}}

\multiput(34,28)(2,4){8}{\circle{2}}
\multiput(34,28)(-2,-4){5}{\circle{2}}

\multiput(38,26)(4,8){4}{\circle{4}}
\multiput(38,26)(-4,-8){2}{\circle{4}}

\multiput(46,22)(8,16){2}{\circle{8}}
\multiput(46,22)(-8,-16){1}{\circle{8}}
\end{picture}
\caption{``Circles''}
\end{center}
\end{figure}

The following theorem is due to  \cite{gilbertnahmod} in the case
$n=3$ and to $\cite{cct}$ in general:

\begin{theorem}\label{multilinear_theorem}
Assume $k:=\dim(\Gamma)<n/2$, and assume that $\Gamma$ is
non-degenerate in the sense that for any $1\le i_1<\dots<i_k\le n$
the space $\Gamma$ is the graph of a function in the variables
$\xi_{i_1},\dots,\xi_{i_k}$. Assume $m$ satisfies
(\ref{new_symbol}). Then $\Lambda$ as in (\ref{multilinear_form})
satisfies (\ref{paraproduct}) whenever $\sum 1/p_j=1$ and
$1<p_j\le \infty$ for all $p_j$.
\end{theorem}

We remark that it is unknown whether the condition
$\dim(\Gamma)<n/2$ can be relaxed in this theorem.

The forms $\Lambda$ have dual multilinear operators. Theorem
\ref{multilinear_theorem} implies a priori estimates for these
multilinear operators. Moreover, these multilinear operators
satisfy estimates which cannot be formulated in terms of $L^p$
estimates for $\Lambda$. Let $(p_1,\dots,p_n)$ be a tuple of real
numbers or $\infty$ such that at most one of these numbers is
negative. If all of them are nonnegative, we say $\Lambda$ is of
type $(p_1,\dots,p_n)$ if (\ref{paraproduct}) holds. If one of
them, say $p_j$, is negative, then we define the dual operator $T$
by
$$\Lambda(f_1,\dots,f_n)=\int T(f_1,\dots, f_{j-1},f_{j+1},\dots, f_n)(x)f_j(x)\, dx. $$
We then say that $\Lambda$ is of type $(p_1,\dots,p_n)$ if
$$\|T(f_1,\dots,f_{j-1},f_{j+1},\dots, f_n)\|_{p_j'}\le C\prod_{i\neq j}\|f_i\|_{p_i}$$ where $p_j'=p_j/(p_j-1)$. Observe $0<p_j'<1$. The following theorem is again due to \cite{gilbertnahmod} ($n=3$) and \cite{cct}:

\begin{theorem}\label{below_one_theorem}
Let $\Gamma$ and $\Lambda$ be as in Theorem
\ref{multilinear_theorem}. Then $\Lambda$ is of type $(p_1,\dots,
p_n)$ if $\sum_j{1/p_j}=1$, at most one of the $p_j$ is negative,
none of the $p_j$ is in $[0,1]$, and
$$1/p_{i_1}+\dots+1/p_{i_r}< \frac{n-2{\rm dim}(\Gamma)+r}2$$
for all $1\le i_1<\dots<i_r\le n$ and $1\le r\le n$.
\end{theorem}

A basic example of a modulation invariant form $\Lambda$ is when
$n=3$ and $m(\xi_1,\xi_2)$ is constant on both sides of a line
$\Gamma$ but not globally constant. With proper choice of
constants this form can be written as
$$\Lambda_\alpha (f_1,f_2,f_3)=\int B_\alpha (f_1,f_2)(x)f_3(x)\, dx$$
with the bilinear Hilbert transform
$$B_\alpha= p.v. \int f_1(x-t)f_2(x-\alpha t)\frac 1t \, dt$$
and a (projective) parameter $\alpha$ determining the direction of
the line $\Gamma$. Theorems \ref{multilinear_theorem} and
\ref{below_one_theorem} in this special case are due to
\cite{laceyt1} and \cite{laceyt2}.

For the bilinear Hilbert transform nondegeneracy specializes to
the condition $\alpha\notin \{0,1,\infty\}$, and the conclusion of
both theorems can be summarized to
\begin{equation}\label{bht}
\|B_\alpha(f_1,f_2)\|_p\le C_{p_1,p_2}\|f_1\|_{p_1}\|f_2\|_{p_2}
\end{equation}
provided $1<p_1,p_2\le \infty$ and $2/3<p<\infty$. The set of
types of such $\Lambda_\alpha$ is the convex hull of the open
triangles $a,b,d$ in Figure ``Hexagon'' which depicts the plane of
$(1/p_1,1/p_2,1/p_3)$ with $\sum_j1/p_j=1$. It is unknown whether
the type-region of $\Lambda_\alpha$ extends to the open triangle
$e$ and its symmetric counterparts.

\setlength{\unitlength}{1.1mm}
\begin{figure}\label{hexagon}
\begin{center}
\begin{picture}(92,73)
\put(46,40){\circle*{0.8}} \put(46,60){\circle*{0.8}}
\put(28,30){\circle*{0.8}} \put(64,30){\circle*{0.8}}
\put(46,30){\circle*{1.6}} \put(19,45){\circle*{0.8}}
\put(28,60){\circle*{0.8}} \put(64,60){\circle*{0.8}}
\put(73,45){\circle*{0.8}} \put(37,15){\circle*{0.8}}
\put(55,15){\circle*{0.8}} \put(37,45){\circle*{1.6}}
\put(55,45){\circle*{1.6}}

\put(5,28){$(1,0,0)$} \put(73,28){$(0,1,0)$}
\put(40,64){$(0,0,1)$}

\put(27,50){$a$} \put(63,50){$b$} \put(45,35){$c$}
\put(45,20){$d$} \put(45,7){$e$}

\multiput(10,60)(6,0){2}{\line(1,0){4}}
\multiput(29,60)(6,0){6}{\line(1,0){4}}
\multiput(82,60)(-6,0){2}{\line(-1,0){4}}
\put(19.4,45){\line(1,0){53.2}} \put(28.4,30){\line(1,0){35.2}}
\multiput(38,15)(6,0){3}{\line(1,0){4}}
\multiput(10,60)(3,-5){2}{\line(3,-5){2}}
\multiput(19.5,44.17)(3,-5){6}{\line(3,-5){2}}
\multiput(46,0)(-3,5){2}{\line(-3,5){2}}
\put(28.2,59.67){\line(3,-5){26.6}}
\put(46.2,59.67){\line(3,-5){17.6}}
\multiput(64.5,59.17)(3,-5){3}{\line(3,-5){2}}
\multiput(19.5,45.83)(3,5){3}{\line(3,5){2}}
\put(28.2,30.33){\line(3,5){17.6}}
\put(37.2,15.33){\line(3,5){26.6}}
\multiput(46,0)(3,5){2}{\line(3,5){2}}
\multiput(55.5,15.83)(3,5){6}{\line(3,5){2}}
\multiput(82,60)(-3,-5){2}{\line(-3,-5){2}}

\put(21,45){\line(3,5){8}} \put(23,45){\line(3,5){7}}
\put(25,45){\line(3,5){6}} \put(27,45){\line(3,5){5}}
\put(29,45){\line(3,5){4}} \put(31,45){\line(3,5){3}}
\put(33,45){\line(3,5){2}} \put(35,45){\line(3,5){1}}

\put(39,45){\line(3,-5){8}} \put(41,45){\line(3,-5){7}}
\put(43,45){\line(3,-5){6}} \put(45,45){\line(3,-5){5}}
\put(47,45){\line(3,-5){4}} \put(49,45){\line(3,-5){3}}
\put(51,45){\line(3,-5){2}} \put(53,45){\line(3,-5){1}}

\put(57,45){\line(3,5){8}} \put(59,45){\line(3,5){7}}
\put(61,45){\line(3,5){6}} \put(63,45){\line(3,5){5}}
\put(65,45){\line(3,5){4}} \put(67,45){\line(3,5){3}}
\put(69,45){\line(3,5){2}} \put(71,45){\line(3,5){1}}

\put(39,15){\line(3,5){8}} \put(41,15){\line(3,5){7}}
\put(43,15){\line(3,5){6}} \put(45,15){\line(3,5){5}}
\put(47,15){\line(3,5){4}} \put(49,15){\line(3,5){3}}
\put(51,15){\line(3,5){2}} \put(53,15){\line(3,5){1}}
\end{picture}
\caption{``Hexagon''}
\end{center}
\end{figure}

We point out a related result by M. Lacey \cite{lacey}:
\begin{theorem}
The maximal truncations of the bilinear Hilbert transform,
\[B_\alpha^{\max}(f,g)(x):=
\sup_{\epsilon>0}\left|\int_{\R\setminus [-\epsilon,\epsilon]}
f(x-t)g(x-\alpha t) \frac 1t \, dt\right|\] also satisfy
(\ref{bht}) provided $\alpha$ is not degenerate.
\end{theorem}

This is stronger than the bounds for the bilinear Hilbert
transform itself.

The main difference in proving the theorems in this section
compared to the discussion in Section \ref{multilinear_section} is
that it is not sufficient to split the functions $f_k$ into
frequency parts supported in $B(0,2^j)\setminus B(0,2^{j-2})$. The
special role that is attributed to the zero frequency by this
splitting is obsolete in the modulation invariant setting. Instead
one has to consider frequency bands of $f_k$ away from the origin
and very narrow, such as intervals $[N-\epsilon, N+\epsilon]$ for
large $N$ and small $\epsilon$. Geometrically these bands can be
viewed as the projections of the circles in Figure ``Circles''
onto the projected coordinate axes. Handling thin frequency bands
requires a new set of techniques. Prior to the work \cite{laceyt1}
and \cite{laceyt2} these techniques have been pioneered in
\cite{carleson} and \cite{fefferman} where the Carleson operator
$$Cf(x)=\sup_\xi |p.v. \int e^{iy\xi}f(x-y)\frac 1{y}\, dy|$$
has been estimated. Note that this operator is modulation
invariant, $C(f)=C(M_\eta f)$. See also \cite{laceyt3}. Most
theorems discussed in this survey have a simpler but significant
model theorem in the dyadic setting, see for example
\cite{mtt:walshbiest}, \cite{thiele''},

\section*{3. Uniform estimates}\setzero \addsec

\vskip-5mm \hspace{5mm}

Theorem \ref{multilinear_theorem} excludes certain degenerate
subspaces $\Gamma$. For some degenerate $\Gamma$ the multilinear
forms split into simpler objects and one can provide $L^p$
estimates also in these degenerate cases; we will give examples
below. This raises the question whether one can prove bounds on
$\Lambda$ uniformly in the choice of $\Gamma$, as $\Gamma$
approaches one of these degenerate cases.

Substantial progress on this question has only been made in the
case $\dim(\Gamma)=1$.
\begin{theorem}\label{uniform_theorem}
Let $n\ge 3$ and $(\eta_1,\dots,\eta_n)$ be a unit vector spanning
the space $\Gamma$, and assume $\eta_j\neq 0$ for all $j$. Define
the metric
$$d(x,y):=\sup_{1\le j\le n}\frac {|x_i-y_i|}{|\eta_i|}$$
and write $d(x,\Gamma):=\inf_{y\in \Gamma} d(x,y)$. Suppose $m$
satisfies the estimate
\begin{equation}\label{ellipse_symbol}
\partial_{\eta'}^\alpha m(\eta')\le \prod_{j=1}^n (\eta_j d(\eta,\Gamma))^{-\alpha_j}
\end{equation}
for all partial derivatives $\partial_{\eta'}^{\alpha}$ up to
order $N$. Then (\ref{paraproduct}) holds for all $2<p_j<\infty$
with $\sum_j 1/p_j=1$ with the bounds uniform in the choice of
$\Gamma$.
\end{theorem}

We discuss uniform estimates for the special case of the bilinear
Hilbert transform. The degenerate directions for $\Gamma$ occur
when the vector $\eta$ is perpendicular to one of the three
projected coordinate axes (see Figure ``Circles''). One of the
degenerate cases ($\alpha=1$) gives rise to the operator
\[B_1(f_1,f_2)=H(f_1\cdot f_2)\]
(Hilbert transform of the pointwise product) or its dual operators
\[ f_2\cdot H(f_3)\ , \ \ f_1\cdot H(f_3). \]
Besides the usual homogeneity $\sum_j1/p_j=1$, the only constraint
for these operators to be of type $(p_1,p_2,p_3)$ is
$1<p_3<\infty$. In Figure ``Hexagon'' this region is the strip
bounded by the horizontal lines through $(0,0,1)$ and $(1,0,0)$.

Thus one expects the constants in the $L^p$ estimates to be
uniform as $\alpha$ approaches $1$ in the intersection of this
strip and the convex hull of triangles $a,b,d$. The above theorem
provides uniform estimates in the inner triangle $c$. This special
case of Theorem \ref{uniform_theorem} was previously shown by
Grafakos/Li \cite{grafakosli}, and Li \cite{li} has shown uniform
estimates in triangles $a$ and $b$. These results together give
uniform bounds in the convex hull of $a,b,c$. Uniform estimates
near the points $(1,0,0)$ and $(0,1,0)$ remain an open question.
Prior to the work of Grafakos/Li \cite{grafakosli}, weak type
uniform bounds were shown \cite{thiele}, \cite{thiele3} in the
common boundary point of triangles $a$ and $c$ (and by symmetry
also $b$ and $c$).

The multiplier condition (\ref{ellipse_symbol}) gives essentially
constant multipliers on regions adapted to the slope of $\Gamma$,
see Figure ``Ellipses''. Observe that all ellipses at a given
scale project essentially onto disjoint regions when projected to
any one one of the coordinate axes. Handling these adapted regions
uniformly requires considerable refinements of the arguments in
\cite{laceyt1} and \cite{laceyt2}.

\setlength{\unitlength}{1.2mm}
\begin{figure}\label{ellipses}
\begin{center}
\begin{picture}(60,62)
\put(30,30){\line(1,0){28}} \put(30,30){\line(-3,5){14}} \put(30,30){\line(-3,-5){14}}

\put(30,30){\line(1,5){6}} \put(30,30){\line(-1,-5){6}}

\put(12,56){$\xi_3$} \put(12,2){$\xi_2$} \put(60,30){$\xi_1$}

\multiput(33,29.4)(1,5){6}{\oval(0.75,4)} \multiput(33,29.4)(-1,-5){6}{\oval(0.75,4)}

\multiput(36,28.8)(2,10){3}{\oval(1.5,8)} \multiput(36,28.8)(-2,-10){3}{\oval(1.5,8)}

\multiput(42,27.6)(4,20){2}{\oval(3,16)} \multiput(42,27.6)(-4,-20){1}{\oval(3,16)}


\end{picture}
\caption{``Ellipses''}
\end{center}
\end{figure}
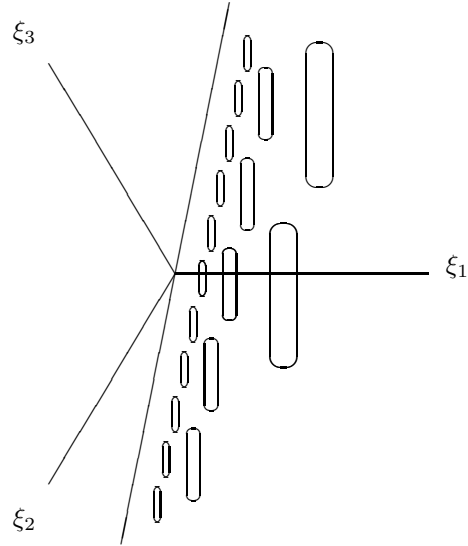

We mention that closely related to the topic of uniform estimates
for the bilinear Hilbert transform is that of bilinear multiplier
estimates for multipliers which are singular along a curve rather
than a line, provided the curve is tangent to a degenerate
direction. Results for such multipliers have been found by Muscalu
\cite{muscalu} and Grafakos/Li \cite{grafakosli2}.

We conclude this section with a remark on the history of the
bilinear Hilbert transform. Calderon is said to have considered
the bilinear Hilbert transform in the 1960's while studying what
has been named Calderon's first commutator. This is the bilinear
operator
\[{\cal C}(A,f)(x)=
p.v. \int \frac {A(x)-A(y)}{(x-y)^{2}}f(y)\,  dy. \]
It can be viewed as a bilinear operator in the derivative
$A'$ of $A$ and the function $f$, and as such has a multplier form as in (\ref{multilinear_form}). To see this, we
can write $C(A,f)$ in terms of $A'$ as a superposition of bilinear Hilbert transforms:
\[{\cal C}(A,f)(x)=p.v. \int \int_{0}^1 A'(x+\alpha(y-x)) \frac 1{x-y}f(y)\,d\alpha dy \]
\[=\int_0^1 B_\alpha(f,A')(x)\, d\alpha. \]
The estimate Calderon was looking for was
\begin{equation}\label{commest}
\|C(A,f)\|_2\le \|A'\|_\infty\|f\|_2.
\end{equation}
Thus he needed good control over the constant $C_\alpha$ as
$\alpha$ approaches $0$ or $1$. However, even finiteness of
$C_\alpha$ was not known to Calderon. Sufficiently good control
over $C_\alpha$ was first established in \cite{thiele}.

The multiplier of ${\cal C}(A',f)$ is more regular than that of
the bilinear Hilbert transform, and Calderon, quitting his
attempts to estimate the bilinear Hilbert transform, proved
estimate (\ref{commest}) by refinements of the methods in Section
\ref{multilinear_section} (see \cite{calderon1}).

\section*{4. More multilinear operators}\setzero \addsec

\vskip-5mm \hspace{5mm}

Theorem \ref{multilinear_theorem} discusses multipliers singular
at a single subspace $\Gamma'$. Cut and paste arguments easily
allow to generalize the theorem to the case of multipliers
singular at finitely many subspaces
${\Gamma_1}',\dots,{\Gamma_k}'$, provided each subspace satisfies
the dimension and non-degeneracy conditions of Theorem
\ref{multilinear_theorem}.

Interesting phenomena occur for multipliers singular at several
subspaces ${\Gamma_1}',$ $\dots,{\Gamma_k}'$ which do not satisfy
the conditions of Theorem \ref{multilinear_theorem}. Some
operators corresponding to multipliers singular at degenerate
subspaces can be written in terms of pointwise products and lower
degree operators and thus can be trivially shown to satisfy $L^p$
estimates. If $m$ is singular at several such subspaces, the
trivial splitting may no longer be possible, and one has to do a
much more subtle analysis.

We consider the special case when the spaces
${\Gamma_1}',\dots,{\Gamma_k}'$ are hyperplanes and the multiplier
is the characteristic function of one of the infinite simplices
been cut out of $\R^n$ by these hyperplanes, see Figure ``Wedge''.
A basic example is the trilinear operator
$$T(f_1,f_2,f_3)(x)=\int _{\alpha_1\xi_1<\alpha_2\xi_2<\alpha_3\xi_3}
\prod_{j=1}^3 \widehat{f}_j(\xi_j)e^{2\pi i x \xi_j}\, d\xi_j$$
and its associated fourlinear form
\begin{equation}\label{biest_lambda}
\Lambda(f_1,f_2,f_3,f_4)=\int_{\sum_{j=1}^4\xi_j=0, \alpha_1\xi_1<\alpha_2\xi_2<\alpha_3\xi_3}\prod_{j=1}^4
\widehat{f}_j(\xi_j)\, d\sigma.
\end{equation}
Here $\alpha_1,\alpha_2,\alpha_3$ are real parameters. If we had
only one of the two constraints $\alpha_1\xi_1<\alpha_2\xi_2$ or
$\alpha_2\xi_2<\alpha_3\xi_3$, then these operators would
decompose trivially.

There is a Zariski open set of values of
$(\alpha_1,\alpha_2,\alpha_3)$ for which $\Lambda$ and $T$ are
well behaved. The following theorem proved in
\cite{mtt:fourierbiest} states such estimates for the generic
point $(1,1,1)$.

\begin{figure}\label{wedge}
\setlength{\unitlength}{1.8mm}
\begin{center}
\begin{picture}(56,30)
\put(20,5){\line(5,1){15}} \multiput(20,5)(-3,3){5}{\line(-1,1){2}}
\multiput(35,8)(-3,3){5}{\line(-1,1){2}} \multiput(30,15)(3,3){2}{\line(1,1){2}}
\multiput(35,8)(3,3){5}{\line(1,1){2}} \multiput(6,19)(4,0.8){4}{\line(5,1){2}}
\multiput(35,20)(4,0.8){4}{\line(5,1){2}} \put(0,10){$\xi_1=\xi_2$} \put(45,13){$\xi_2=\xi_3$}
\end{picture}
\caption{``Wedge''}
\end{center}
\end{figure}

\begin{theorem}\label{generic}
For $\alpha_1,\alpha_2,\alpha_3=1$ the form $\Lambda$ as in
(\ref{biest_lambda}) satisfies estimates
$$\Lambda(f_1,f_2,f_3,f_4)\le C_{p_1,\dots,p_4}\prod_{j=1}^4\|f_j\|_{p_j}$$
if $1< p_j<\infty$ and $\sum_j1/p_j=1$. The trilinear form $T$
satisfies in addition estimates mapping into $L^p$ with $p<1$, in
particular
$$\|T(f_1,f_2,f_3)\|_{2/3}\le C \prod_{j=1}^3 \|f_j\|_2. $$
\end{theorem}

An example for a degenerate choice of
$(\alpha_1,\alpha_2,\alpha_3)$ is $(1,-1,1)$. In this case there
is a negative result \cite{mtt:counterexample}:

\begin{theorem}\label{counterexample}
For $\alpha_1=1,\alpha_2=-1,\alpha_3=1$ the a priori estimate
$$\|T(f_1,f_2,f_3)\|_{2/3}\le C \prod_{j=1}^3 \|f_j\|_2$$
does not hold.
\end{theorem}

Theorem \ref{counterexample} is proved by applying $T$ to
functions $f_1,f_2,f_3$ which are suitable truncations of
imaginary Gaussians (chirps) $e^{i\beta x^2}$. The operator of
Theorem \ref{counterexample} appears naturally in eigenfunction
expansions of one dimensional Schr\"odinger operators, see the
work of Christ/Kiselev \cite{ck},\cite{ck-2}. A positive result on
discrete models of these expansions using the modulation invariant
theory can be found in \cite{mtt:walshnlc}.

\label{lastpage}

\end{document}